\documentclass{article}
\usepackage[T1]{fontenc}
\usepackage[latin1]{inputenc}
\usepackage{amssymb}

\makeatletter

\providecommand{\LyX}{L\kern-.1667em\lower.25em\hbox{Y}\kern-.125emX\@}

\usepackage{amsmath}
\usepackage{amsthm}
 \newtheorem{thm}{Theorem}[section]

 \usepackage{verbatim}
 \theoremstyle{definition}
 \newtheorem{defn}[thm]{Definition}
 \theoremstyle{definition}
  \newtheorem{example}[thm]{Example}
 \theoremstyle{plain}    
 \newtheorem{fact}[thm]{Fact} 
 \theoremstyle{plain}    
 \newtheorem{lem}[thm]{Lemma} 
 \theoremstyle{remark}    
 \newtheorem{note}[thm]{Note} 
 \newcommand{\lyxrightaddress}[1]{
   \par {\raggedleft \begin{tabular}{l}\ignorespaces
   #1
   \end{tabular}
   \vspace{1.4em}
   \par}
 }

\date{\today}
\usepackage[english]{babel}
\makeatother
\begin{document}

\newcommand{\norm}[1]{\lVert #1\rVert }

\newcommand{\abs}[1]{\lvert #1\rvert }

\newcommand{\LL}{\langle \! \langle }

\newcommand{\GG}{\rangle \! \rangle }

\newcommand{\bform}[2]{\left\langle #1\left|\right.#2\right\rangle }

\newcommand{\bbform}[2]{\left\langle \! \left\langle #1\left|\right.#2\right\rangle \! \right\rangle }

\newcommand{\emu}{e_{\mu }(\xi ;\cdot )}

\newcommand{\enu}{e_{\nu }(\eta ;\cdot )}

\newcommand{\epq}{E_{p,q}^{1}}

\newcommand{\epa}{[E_{p}]_{\alpha }}

\newcommand{\fnp}{\left|f_{n}\right|_{p}^{2}}

\newcommand{\holo}{\mathrm{Hol}_{0}(\mathcal{N}_{\mathbb{C}})}

\newcommand{\holonn}{\mathrm{Hol}_{0}(\mathcal{N}_{\mathbb{C}},\mathcal{N}_{\mathbb{C}})}

\newcommand{\hpq}{(H_{p,q,\mu })}

\newcommand{\hpqnu}{(H_{p,q,\nu })}

\newcommand{\hpqo}{(H_{p_{o},q_{o},\mu })}

\newcommand{\Hpc}{H_{p,\mathbb{C}}}

\newcommand{\Hmpc}{H_{-p,\mathbb{C}}}

\newcommand{\hmpq}{(H_{-p,-q,\mu })}

\newcommand{\hmpqnu}{(H_{-p,-q,\nu })}

\newcommand{\hmpqo}{(H_{-p_{o},-q_{o},\mu })}

\newcommand{\Ncp}{\mathcal{N}_{\mathbb{C}}^{\prime }}

\newcommand{\Ncpt}[1]{\mathcal{N}_{\mathbb{C}}^{\prime \hat{\otimes }#1}}

\newcommand{\Nct}[1]{\mathcal{N}_{\mathbb{C}}^{\hat{\otimes }#1}}

\newcommand{\test}{(\mathcal{N})^{1}}

\newcommand{\dist}{(\mathcal{N})^{-1}}

\newcommand{\Nc}{\mathcal{N}_{\mathbb{C}}}

\newcommand{\Pn}{P_{n,\mu }}

\newcommand{\Pm}{P_{m,\mu }}

\newcommand{\Pna}{P_{n,\mu ,\alpha }}

\newcommand{\Qn}{Q_{n,\mu }}

\newcommand{\Qm}{Q_{m,\nu }}

\newcommand{\Qna}{Q_{n,\mu ,\alpha }}

\newcommand{\Qma}{Q_{m,\nu ,\alpha }}

\title{OPERATORS ON SPACES GENERATED BY INFINITE DIMENSIONAL APPELL POLYNOMIALS }

\author{EUGENE YABLONSKY}

\maketitle

\lyxrightaddress{Department of Mathematical Sciences \\
Worcester Polytechnic Institute, 100 Institute Road\\
Worcester, MA 01609, USA\\
eugyabl@wpi.edu\\
}

\begin{abstract}
It is known that many constructions arising in the classical Gaussian
infinite dimensional analysis can be extended to the case of more
general measures. One such extension can be obtained through biorthogonal
systems of Appell polynomials and generalized functions. In this paper,
we consider linear continuous operators from a nuclear Fréchet space
of test functions to itself in this more general setting. We construct
an isometric integral transform (biorthogonal CS-transform) of those
operators into the space of  germs of holomorphic functions on a locally
convex infinite dimensional nuclear space. Using such transform, we
provide characterization theorems and give biorthogonal chaos expansion
for operators.
\end{abstract}
\setcounter{tocdepth}{2}

\setcounter{secnumdepth}{3}

\section*{Introduction}

Recently, many researchers\cite{hida:2002,aldakos96,akk2001,kosis97}
have been working on extending Gaussian Infinite Dimensional Analysis
and White Noise Calculus beyond the case of the Gaussian measure.
One possible approach to this problem is through \emph{biorthogonal}
systems of polynomials and generalized functions. That approach was
discussed by Yu.~Daletsky, S.~Albeverio, Yu.~Kondratiev, L.~Streit,
W.~Westerkamp, J.-A.~Yan, J.~Silva, et al.\cite{koswy95,kosis97,dal:91}

To illustrate their idea, consider a measure $d\mu (x)=p(x)dx$ on
$\mathbb{R}$, where $p(x)$ is a smooth positive function of class
$L_{1}(\mathbb{R},\exp (\varepsilon |x|)dx)$, $\varepsilon >0$.
In case of the Gaussian measure, the Taylor expansion of the normalized
exponential function \textbf{$e_{\mu }(\xi ;x):=\frac{\exp (x\xi )}{\int _{\mathbb{R}}\exp (x\xi )d\mu (x)}=\sum _{n=0}^{\infty }\frac{\xi ^{n}}{n!}P_{n}(x)$}
generates an orthogonal system of Hermite polynomials. In case of
a more general measure, the Appell polynomials $P_{n}$ are not necessarily
orthogonal; however, one can construct a dual system $Q_{n}(x)=\left(\frac{d^{n}}{dx^{n}}\right)^{*}\mathbf{1}_{\mathbb{R}}=(-1)^{n}\frac{p^{(n)}(x)}{p(x)}$.
This biorthogonal system was introduced for smooth measures by Dalecky
et al.\cite{dal:91,aldakos96} and later extended to a broader class
of non-degenerate measures with analytic characteristic functionals
by Kondratiev et al.\cite{koswy95,kosis97}

Using the biorthogonal system, the mentioned authors constructed the
nuclear spaces $\left(\mathcal{N}\right)^{\beta }$ of {}``Appell''
test functions and the dual spaces $\left(\mathcal{N}\right)^{-\beta }$,
$0\leq \beta \leq 1$, with the space $\dist $ being the largest
and the most technically challenging. All those spaces were characterized
in terms of so called $C_{\mu }$- and $S_{\mu }$- transforms, which
are biorthogonal analogues of the S-transform in white noise calculus. 

In this paper, we characterize linear continuous operators ${\test \to \test }$.
For any such operator $B$ we define a {}``local $CS_{\mu ,\nu }$-symbol''
${\check{B}_{\mu ,\nu }=C_{\mu }(Be_{\nu })}$, where $\mu ,\nu $
are two measures on a conuclear space $\mathcal{N}'$, and $e_{\nu }$
is a normalized exponent. The definition of that symbol needs to be
carefully interpreted, as $e_{\nu }$ and $Be_{\nu }$ belong not
to $\test $ but to larger Hilbert spaces. The symbol $\check{B}_{\mu ,\nu }$
is actually a germ of complex-valued functions holomorphic in \emph{cylindrical}
0-neighborhoods of $\Nc \times \Nc $. We show that under certain
conditions it uniquely characterizes the operator $B$. Our approach
is similar in spirit to one used for characterization of white noise
operators (see, e.g., Obata.\cite{obata94}), however our operator
symbols are holomorphic only locally, in very special cylindrical
domains. 

The case of operators ${\test \to \dist }$ is studied in the author's
thesis\cite{yabl:2003} by analyzing so called local $S_{\mu \nu }$-symbols
${\hat{B}_{\mu ,\nu }=S_{\mu }(Be_{\nu })}$.

The paper is organizes as follows. Section 1 recalls definitions related
to locally holomorphic and entire functions on nuclear spaces. Section
2 describes the spaces of test and generalized functions introduced
by Kondratiev et. al.\cite{koswy95,aldakos96,kosis97} Those spaces
are defined by mean of biorthogonal Appell systems. Section 3 introduces
local symbols of operators and gives a characterization for operators
$\test \to \test $.

\section{Preliminaries}

\subsection{Nuclear spaces}

Let $\mathcal{N}\subset H\subset \mathcal{N}'$ be a real nuclear
Gel'fand triple and $\langle \cdot |\cdot \rangle $ the canonical
bilinear form on $\mathcal{N}\times \mathcal{N}'$ such that $\langle \cdot |\cdot \rangle $
is an extension of the inner product $(\cdot ,\cdot )$ on $H$. Denoted
by $\left|\cdot \right|$ the norm on $H$. Without loss of generality,
we can always assume that the nuclear Fréchet space $\mathcal{N}$
is represented as a projective limit of a sequence ${H=H_{0}\supset H_{1}\supset H_{2}\supset ...}$
of Hilbert spaces with norms ${|\cdot |=|\cdot |_{0}\leq |\cdot |_{1}\leq |\cdot |_{2}\leq ...}$;
moreover, for every $H_{p}$ there is $p'>p$ such that the embedding
${I_{p',p}:H_{p'}\hookrightarrow H_{p}}$ is a Hilbert-Schmidt operator. 

Let $H_{-p}:=H_{p}^{\prime }$ be the dual to the Hilbert space $H_{p}$.
Then the space dual to $\mathcal{N}$ is given by ${\mathcal{N}'=\textrm{ind}\lim _{p\to \infty }H_{-p}}$.
See Schaefer\cite{schaefer71} for details. We will denote the Hilbert
norm on $H_{-p}$ by $\abs{\cdot }_{-p}$.

Tensor product of nuclear spaces\index{Tensor product of nuclear spaces}
(or $\pi $-product) $\mathcal{N}\otimes \mathcal{M}$ is the nuclear
space defined as the projective limit of Hilbert space tensor products
${H_{p}\otimes L_{p}}$, $p=0,1,2,\ldots $. Similarly, we define
tensor powers $\mathcal{N}^{\otimes n}$. The space dual to $\mathcal{N}^{\otimes n}$
can be presented as the inductive limit ${\mathcal{N}^{\prime \otimes n}:=\textrm{ind}\lim H_{-p}^{\otimes n}}$.
We will preserve the notations $\abs{\cdot }_{p}$ and $\abs{\cdot }_{-p}$
for the norms on the tensor powers $H_{p}^{\otimes n}$ and $H_{-p}^{\otimes n}$
respectively, as well as their complexifications $H_{p,\mathbb{C}}^{\otimes n}$,
and $H_{-p,\mathbb{C}}^{\otimes n}$. Symmetric tensor product will
be denoted by $\hat{\otimes }$.

\subsection{Locally Holomorphic and Entire Functions}

Let $\Nc $ is the complexification of $\mathcal{N}$. We will consider
functions holomorphic in a neighborhood of $0\in \Nc $.

\begin{fact}
\label{cor:locally-holomorphic}A function $G:\Nc \rightarrow \mathbb{C}$
is holomorhic at $0$ if and only if 
\end{fact}
\begin{enumerate}
\item There exist $p$ and $\rho >0$ such that for all $\xi _{0}\in \mathcal{N}_{\mathbb{C}}$
with $\left|\xi _{0}\right|_{p}\leq \rho $ and for all $\xi \in \mathcal{N}_{\mathbb{C}}$
the function of one complex variable ${\lambda \mapsto G(\xi _{0}+\lambda \xi )}$
is analytic at $0\in \mathbb{C}$, and 
\item There exists $c>0$ such that for all $\xi \in \mathcal{N}_{\mathbb{C}}$
with $\left|\xi \right|_{p}\leq \rho $, we have $\left|G(\xi )\right|\leq c$. 
\end{enumerate}
We do not discern between different restrictions of one function,
i.e., we identify $F$ and $G$ if there exists an open 0-neighborhood
$U\subset \mathcal{N}_{\mathbb{C}}$ such that $F(\xi )=G(\xi )$
for all $\xi \in U$. Denote by $\mathrm{Hol}_{0}(\mathcal{N}_{\mathbb{C}})$
the algebra of germs of functions $\mathcal{N}_{\mathbb{C}}\rightarrow \mathbb{C}$
holomorphic at $0$. That algebra is equipped with the inductive limit
topology given by the family of norms\begin{eqnarray*}
\mathsf{n}_{p,l,\infty }(G) & = & \sup _{\abs{\xi }_{p}\leq 2^{-l}}\abs{G(\xi )}\, ,\, \, \, \, \, p,l\in \mathbb{N}.
\end{eqnarray*}

\begin{defn}
\label{def:Entire}A function $G:\Ncp \to \mathbb{C}$ is an \emph{entire
function of growth ${k\in [1,2]}$ and a minimal type} if $G$ is
holomorphic on every $\Hmpc $, $p=0,1,\ldots $, and for any $\varepsilon >0$
there exists $C>0$ such that ${\abs{G(z)}\leq C\exp (\varepsilon \abs{z}_{-p}^{k})}$,
${z\in \Hmpc }$. 
\end{defn}
Denote by $\mathcal{E}_{\min }^{k}(\Ncp )$ the space of all such
functions. The space is endowed with the projective limit topology
with respect to the countable system of norms\begin{eqnarray*}
\mathsf{n}_{p,l,k}(G) & = & \sup _{z\in H_{-p,\mathbb{C}}}\abs{G(z)}\exp (-2^{-l}\abs{z}_{-p}^{k}),\, \, \, p,l\in \mathbb{N}.
\end{eqnarray*}
For more details about these spaces see, e.g., Kondratiev et al.\cite{koswy95}
and the textbook of Dineen\cite{dineen81}.

\section{Biorthogonal Appell System}

\begin{defn}
A function $\phi :\mathcal{N}'\rightarrow \mathbb{C}$ of the form
\[
\phi (x)=\sum _{n=0}^{N}\langle x^{\otimes n}|\phi _{n}\rangle \, \, \, \, \, \phi _{n}\in \mathcal{N}_{\mathbb{C}}^{\widehat{\otimes }n}\]
 is called a \emph{continuous polynomial} on $\mathcal{N}'$. The
set of all continuous polynomials is denoted by $\mathcal{P}(\mathcal{N}')$,
while the dual space is denoted by \textbf{}$\mathcal{P}'(\mathcal{N}')$.
\end{defn}
Consider the $\sigma $-algebra $C(\mathcal{N}')$ generated by cylindrical
sets on $\mathcal{N}'$, which coincides with the Borel $\sigma $-algebras
generated by the strong and inductive limit topology. We consider
the class of measures $\mu $ on $C(\mathcal{N}')$ satisfying to
the following assumptions:

\begin{itemize}
\item ASSUMPTION 1. The measure $\mu $ on the space $\mathcal{N}'$ has
an analytic Laplace transform in a neighborhood of $0\in \mathcal{N}_{\mathbb{C}}$\[
\textrm{E}_{\mu }(\exp \langle \cdot |\xi \rangle )=\int _{\mathcal{N}'}e^{\langle x|\theta \rangle }d\mu (x)\, \in \, \mathrm{Hol}_{0}(\mathcal{N}_{\mathbb{C}});\]

\item ASSUMPTION 2. The measure $\mu $ is non-degenerate, that is for every
continuous polynomial $\phi \in \mathcal{P}(\mathcal{N}')$\begin{eqnarray*}
\phi =0\, \, \, \mu -\textrm{almost everywhere} & \Longleftrightarrow  & \textrm{ }\phi \equiv 0.
\end{eqnarray*}

\end{itemize}
It follows from the Assumption 1 that $\mathcal{P}(\mathcal{N}')$
is densely embedded into $L^{2}(\mu )$ (see the book of Skorohod\cite{Sk74},
Section 10 for details). We obtain the Gel'fand triple $\mathcal{P}(\mathcal{N}')\subset L^{2}(\mu )\subset \mathcal{P}'(\mathcal{N}')$.
The bilinear dual pairing between $\mathcal{P}(\mathcal{N}')$ and
$\mathcal{P}'(\mathcal{N}')$ with respect to $\mu $ is denoted by
$\bbform{\cdot }{\cdot }_{\mu }$; for $\phi \in L^{2}(\mu )$ and
$\psi \in \mathcal{P}(\mathcal{N}')$ we have $\bbform{\phi }{\psi }_{\mu }=\int _{\mathcal{N}'}\phi (x)\psi (x)d\mu (x)$.

\subsection{\label{sub:Appell-Polynomials}Normalized Exponentials and Appell
Polynomials}

Consider a normalized exponent\begin{eqnarray*}
e_{\mu }(\xi ;x) & := & \frac{\exp \langle x|\xi \rangle }{\textrm{E}_{\mu }(\exp \langle x|\xi \rangle )}\, ,\, \, \, \, \xi \in \mathcal{N}_{\mathbb{C}},x\in \mathcal{N}_{\mathbb{C}}^{\prime }\, .
\end{eqnarray*}
 It is a well defined function for $\xi $ in some $0$-neighborhood
$U_{0}\subset \mathcal{N}_{\mathbb{C}}$.

For every $x\in \Ncp $, the function $\xi \mapsto e_{\mu }(\xi ;x)$
is holomorphic (and also bounded) on some 0-neighborhood $U_{0}=\{\xi \in \mathcal{N}_{\mathbb{C}};\left|\xi \right|_{p_{0}}<2^{-q_{0}}\}.$
Consider its Taylor series about the point $\xi =0$\begin{eqnarray*}
e_{\mu }(\xi ;x) & = & \sum _{n=0}^{\infty }\frac{1}{n!}d^{n}e_{\mu }(0;x)\xi ^{n},
\end{eqnarray*}
where $d^{n}e_{\mu }(0;x)$ is a symmetric n-form, and we use the
notation\[
d^{n}e_{\mu }(0;x)\xi ^{n}:=d^{n}e_{\mu }(0;x)(\xi ,\ldots ,\xi )\, .\]
Using the polarization identity, we obtain that the form ${d^{n}e_{\mu }(0;x)}$
is $H_{p,\mathbb{C}}$-continuous on the whole ${\mathcal{N}_{\mathbb{C}}^{\hat{\otimes }n}}$,
i.e.\[
\abs{\frac{1}{n!}d^{n}e_{\mu }(0;x)(\xi _{1},\ldots ,\xi _{n})}\leq \textrm{const}\prod _{j=1}^{n}\left|\xi _{j}\right|_{p}\, .\]
By the kernel theorem, for all $p'>p$ such that the embedding ${i_{p',p}:H_{p'}\rightarrow H_{p}}$
is a Hilbert-Schmidt operator, there exist unique kernels ${\Pn (x)\in H_{-p',\mathbf{C}}^{\hat{\otimes }n}}$
($x$ is fixed) such that \[
d^{n}e_{\mu }(0;x)\xi ^{n}:=\langle \Pn (x)|\xi ^{\otimes n}\rangle \, .\]
The system $\mathbb{P}_{\mu }=\left\{ \bform{\Pn (\cdot )}{\varphi _{n}};\, \varphi _{n}\in \Nct{n},n=0,1,2,\ldots \right\} $
is called \emph{The Appell System of polynomials} (on $\Ncp )$.

\subsection{\label{sub:Hida-Differential-Operator}Differential operators and
$\mathbb{Q}_{\mu }$-System}

Let $\Phi _{n}\in \Ncpt{n}$; define the differential operator ${D(\Phi _{n}):\mathcal{P}(\mathcal{N}')\rightarrow \mathcal{P}(\mathcal{N}')}$,
acting on monomials $\langle x^{\otimes m}|\phi _{m}\rangle $, $\phi _{m}\in \Nct{m}$
as

\begin{eqnarray*}
D(\Phi _{n})\langle x^{\otimes m}|\phi _{m}\rangle  & := & \frac{m!}{(m-n)!}\langle x^{\otimes (m-n)}\hat{\otimes }\Phi _{n}|\phi _{m}\rangle 
\end{eqnarray*}
 whenever $m>n$ and $0$ otherwise. For $\Phi _{1}\in \mathcal{N}_{\mathbb{C}}'$,
$\varphi \in \mathcal{P}(\mathcal{N}')$ we have ${D(\Phi _{1})\varphi =\left.\frac{d}{dt}\varphi (\cdot +t\Phi _{1})\right|_{t=0}}$.

\begin{lem}
(Kondratiev et al.\cite{koswy95}) $D(\Phi _{n}):\mathcal{P}(\mathcal{N}')\rightarrow \mathcal{P}(\mathcal{N}')$
is a continuous linear operator which acts on on monomials$\langle \Pm |\phi _{m}\rangle $,
$\phi _{m}\in \Nct{m}$ as follows: \[
D(\Phi _{n})\langle \Pm |\phi _{m}\rangle =\frac{m!}{(m-n)!}\langle P_{m-n,\mu }\hat{\otimes }\Phi _{n}|\phi _{m}\rangle \, ,\, \, \, \]
 whenever $m>n$ and $0$ otherwise.\hfill{}$\square $
\end{lem}
For any $\Phi _{n}\in \mathcal{N}_{\mathbb{C}}^{\prime \widehat{\otimes }n}$
define the generalized function $\Qn (\Phi _{n})\in \mathcal{P}'(\mathcal{N}')$

\[
\Qn (\Phi _{n}):=D(\Phi _{n})^{*}\mathbf{1}_{\mathcal{N}'}.\]

\begin{thm}
(Kondratiev et al.\cite{koswy95}) There is biorthogonality relation
w.r.t. $\mu $\[
\bbform{\Qn (\Phi _{n})}{\langle \Pm |\phi _{m}\rangle \, }=\delta _{m,n}n!\langle \Phi _{n}|\phi _{n}\rangle \]
 and every element $\Phi \in \mathcal{P}'(\mathcal{N}')$ can be presented
as ${\Phi =\sum _{n=0}^{\infty }\Qn (\Phi _{n})}$.\hfill{}$\square $
\end{thm}

\subsection{The Test Space of Appell functions $(\mathcal{N})^{1}$}

For each continuous polynomial $\phi (x)=\sum _{n=0}^{N}\langle \Pn (x)|\phi _{n}\rangle \in \mathcal{P}(\mathcal{N}')$,
with $\phi _{n}\in \mathcal{N}_{\mathbb{C}}^{\widehat{\otimes }n}$,
define the norm\[
\left\Vert \phi \right\Vert _{p,q,\mu }^{2}:=\sum _{n=0}^{N}(n!)^{2}2^{qn}\left|\phi _{n}\right|_{p}^{2}\]
Let the Hilbert space $(H_{p,q,\mu })$ be the completion of $\mathcal{P}(\mathcal{N}')$
with respect to the norm $\left\Vert \phi \right\Vert _{p,q,\mu }$.
Define \emph{the test space of Appell functions}\begin{eqnarray*}
(\mathcal{N})^{1}:=\textrm{proj}\lim _{p,q\in \mathbf{N}}(H_{p,q,\mu }) & . & 
\end{eqnarray*}

\begin{thm}
\label{thm:test-space}(Kondratiev et al.\cite{koswy95}) The space
$(\mathcal{N})^{1}$ is a nuclear space densely embedded in $L^{2}(\mu )$.
Every element $\varphi \in \test $ has a unique extension to $\Nc $
as an element of the space of entire functions of growth $1$ and
a minimal type $\mathcal{E}_{\min }^{1}(\Ncp )$ (see the definition
\ref{def:Entire}), and there is the topological identity\[
\test =\mathcal{E}_{\min }^{1}(\mathcal{N}'):=\left\{ \varphi \big |_{\mathcal{N}'};\varphi \in \mathcal{E}_{\min }^{1}(\Ncp )\right\} ,\]
so that $\test $ does not depend on the choice of the measure $\mu $.\hfill{}$\square $
\end{thm}
\begin{example}
Consider the $\mu $-exponential
\end{example}
\[
e_{\mu }(\xi ;\cdot )=\frac{\exp \langle x|\xi \rangle }{\textrm{E}_{\mu }(\exp \langle x|\xi \rangle )}=\sum _{n=0}^{\infty }\frac{1}{n!}\langle \Pn (\cdot )|\xi ^{\otimes n}\rangle \, ,\, \, \, \, \xi \in \mathcal{N}_{\mathbb{C}}\, .\]
Its norm\[
\left\Vert e_{\mu }(\theta ;\cdot )\right\Vert _{p,q,\mu }^{2}=\sum _{n=0}^{\infty }(n!)^{2}2^{nq}\left(\frac{\left|\theta \right|_{p}^{n}}{n!}\right)^{2}\]
 is finite if and only if $\left|\theta \right|_{p}^{2}<2^{-q}$,
so that $e_{\mu }(\theta ;\cdot )\not \in (\mathcal{N}_{\mathbb{C}})^{1}$
whenever $\theta \not =0$. 

However, for any $\theta \in \mathcal{N}_{\mathbb{C}}$, $\left|\theta \right|_{p}^{2}<2^{-q}$
we have $e_{\mu }(\theta ;\cdot )\in (H_{p,q,\mu })$, that is $e_{\mu }(\theta ;\cdot )$
is a test function of finite order.

\begin{lem}
\label{lem:total-set} (Kondratiev et al.\cite{koswy95}) The family
of normalized exponentials ${\left\{ e_{\mu }(\theta ;\cdot );\left|\theta \right|_{p}^{2}<2^{-q},\theta \in \mathcal{N}_{\mathbb{C}}\right\} }$
is a total set in $(H_{p,q,\mu })$. \hfill{}$\square $
\end{lem}

\subsection{Space of Generalized Functions $\dist $}

The Hilbert space $\hmpq :=\hpq ^{*}$ can be presented as a subspace
of $\mathcal{P}'(\mathcal{N}')$ consisting of those elements ${\Phi =\sum _{n=0}^{\infty }\Qn (\Phi _{n})\in \mathcal{P}'(\mathcal{N}')}$
for which all the norms ${\left\Vert \Phi \right\Vert _{-p,-q,\mu }^{2}:=\sum _{n=0}^{\infty }2^{-qn}\left|\Phi _{n}\right|_{-p}^{2}}$
are finite. Define the inductive limit

\[
(\mathcal{N})_{\mu }^{-1}:=\textrm{ind}\lim _{p,q\in \mathbf{N}}\hmpq .\]
The space $(\mathcal{N})_{\mu }^{-1}$ is the dual space of $\test $
relative to $L^{2}(\mu )$. As the inductive limit topology is equivalent
to the strong topology on $(\mathcal{N})_{\mu }^{-1}$, we can drop
the subscript $\mu $, unless we want to discuss a particular representation
of the space $(\mathcal{N})_{\mu }^{-1}$. 

We obtain the nuclear triple $(\mathcal{N})^{1}\subset (L^{2})=L^{2}(\mathcal{N}',\mu )\subset (\mathcal{N})_{\mu }^{-1}$
and have the following series of continuous embeddings:

\[
(\mathcal{N})^{1}\subset ...\subset \hpq \subset (L^{2})_{\mu }\subset \ldots \subset \hmpq \subset ...\subset (\mathcal{N})_{\mu }^{-1}\]

\begin{example}
(Generalized Radon-Nikodim Derivative). For any $\xi \in \mathcal{N}'$
we can define the generalized function\begin{eqnarray*}
\rho _{\mu }(\xi ) & = & \sum _{n=0}^{\infty }(-1)^{n}\Qn \left(\frac{\xi ^{\otimes n}}{n!}\right)\subset \dist \, .
\end{eqnarray*}
It follows from Lemma \ref{lem:total-set} and the following identity
\[
\int _{\mathcal{N}'}e_{\mu }(\eta ;x-\xi )d\mu (x)=\frac{\int _{\mathcal{N}'}e^{\langle x-\xi |\eta \rangle }d\mu (x)}{\int _{\mathcal{N}'}e^{\langle x|\eta \rangle }d\mu (x)}=e^{-\bform{\xi }{\eta }}=\bbform{\rho _{\mu }(\xi )}{e_{\mu }(\eta ;\cdot )}_{\mu }\]
that $\rho _{\mu }(\xi )$ coincides with the Radon-Nikodim derivative
$\frac{d\mu (x+\xi )}{d\mu (x)}$ whenever that derivative exists. 
\end{example}

\subsection{Local $S_{\mu }$-transform of Generalized Functions}

Spaces of test and generalized functions are often characterized by
mean of integral transforms. Following Kondratiev et al.\cite{koswy95},
let us introduce so called $S_{\mu }$-transform (which can be regarded
as a normalized Laplace transform)

\begin{eqnarray*}
\left(S_{\mu }\Phi \right)(\theta ) & := & \bbform{\Phi }{e_{\mu }(\theta ;\cdot )}_{\mu }\, ,\, \, \, \Phi \in \dist \, .
\end{eqnarray*}
This definition needs to be properly interpreted, as $e_{\mu }(\theta ;\cdot )\not \in (\mathcal{N}_{\mathbb{C}})^{1}$
whenever $\theta \not =0$. However, it is possible to define the
$S_{\mu }$-transform locally as follows. 

As every generalized function $\Phi \in \dist $ is of finite order,
there exist $p,\, q>0$ such that ${\Phi \in \hmpq }$. Consider a
0-neighborhood ${U_{p,q}=\{\theta \in \mathcal{N}_{\mathbb{C}};\, \left|\theta \right|_{p}^{2}<2^{-q}\}}$.
For $\theta \in U_{p,q}$ we have $e_{\mu }(\theta ;\cdot )\in (H_{p,q,\mu })$,
as\[
\left\Vert e_{\mu }(\theta ;\cdot )\right\Vert _{p,q,\mu }^{2}=\sum _{n=0}^{\infty }(n!)^{2}2^{nq}\left(\frac{\left|\theta \right|_{p}^{n}}{n!}\right)^{2}<\infty \, ,\]
 Therefore, the bilinear form $\bbform{\Phi }{e_{\mu }(\theta ;\cdot )}_{\mu }$
is well defined on $U_{p,q}$. 

\begin{thm}
(Kondratiev et al.\cite{koswy95}) The $S_{\mu }$-transform is a
topological isomorphism from $\dist $ to $\holo $.\hfill{}$\square $
\end{thm}

\subsection{The $C_{\mu }$-transform of Test Functions}

\begin{defn}
The $C_{\mu }$-transform for a test function $\varphi \in \test $
is defined as\begin{eqnarray*}
\left(C_{\mu }\varphi \right)(\xi ) & :=\int _{\mathcal{N}'}\varphi (x+\xi )d\mu (x)= & \bbform{\rho _{\mu }(-\xi )}{\varphi }_{\mu }\, ,\, \, \, \xi \in \Ncp \, .
\end{eqnarray*}

\end{defn}
\begin{thm}
(Kondratiev et al.\cite{koswy95})\label{thm:c-transfrom} The $C_{\mu }$-transform
is a topological isomorphism from $\test $ to $\mathcal{E}_{\min }^{1}(\Ncp )$.\hfill{}$\square $
\end{thm}
\begin{note}
If $\mu $ is the Gaussian measure, then the $C_{\mu }$- and $S_{\mu }$-
transforms coincide.
\end{note}

\section{Characterization of Operators $\test \to \test $}

Let $B:\test \to \test $ be a continuous linear operator, and $\mu ,\nu $
be non-degenerate measures having analytic Laplace transform. We want
to represent $B$ by a holomorphic function, that is to define a symbol
of $B$.

\subsection{\label{sub:cs-transform}Local $CS_{\mu \nu }$-Symbols}

As $\test $ is continuously embedded into $\dist $, one can consider
the local $S_{\mu \nu }$-symbol ${\hat{B}_{\mu ,\nu }(\xi ,\eta ):=S_{\mu }[Be_{\nu }(\eta ;\cdot )](\xi )}$,
see the author's thesis\cite{yabl:2003}. However, because of biorthogonality
it is more convenient to work with a somewhat different symbol ${\check{B}_{\mu ,\nu }(\xi ,\eta ):=C_{\mu }[Be_{\nu }(\eta ;\cdot )](\xi )=\bbform{\rho _{\mu }(-\xi )}{Be_{\nu }(\eta ;\cdot )}_{\mu }}$,
which we are going to introduce now. That definition needs to be carefully
interpreted, as $e_{\nu }(\eta ;\cdot )$ is not in $\test $ but
rather in one of the approximating Hilbert spaces. The same thing
might happen with $Be_{\nu }(\eta ;\cdot )$, so that we also need
to extend $C_{\mu }$-transform from $\test $ to larger spaces.

The space $\test $ is continuously embedded into every Hilbert space
$\hpqo $, $p_{o},q_{o}\geq 0$. Applying the kernel theorem to the
continuous operator \[
B:\test \to \test \hookrightarrow \hpqo \]
we obtain that there exist integers $p\geq 0,\, q\geq 0$ such that
$B$ is a continuous linear operator between Hilbert spaces ${\hpqnu \rightarrow \hpqo }$. 

Consider a 0-neighborhood $U_{p,q}=\{\eta \in \mathcal{N}_{\mathbb{C}};\, \left|\eta \right|_{p}^{2}<2^{-q-1}\}$.
For any ${\eta \in U_{p,q}}$\[
\left\Vert e_{\nu }(\eta ;\cdot )\right\Vert _{p,q,\nu }^{2}=\sum _{n=0}^{\infty }(n!)^{2}2^{nq}\left(\frac{\left|\eta ^{\otimes n}\right|_{p}}{n!}\right)^{2}=\frac{1}{1-2^{q}\left|\eta \right|_{p}^{2}}<\infty \, ,\]
 so that $e_{\nu }(\eta ;\cdot )\in (H_{p,q,\nu })$. For any $\xi \in \Nc \subset H_{-p_{o},\mathbb{C}}$
we have\begin{eqnarray}
\norm{\rho _{\mu }(-\xi )}_{-p_{o},-q_{o},\mu } & = & \sqrt{\sum _{n=0}^{\infty }2^{-nq_{o}}\left(\frac{\left|\xi ^{\otimes n}\right|_{-p_{o}}}{n!}\right)^{2}}\nonumber \\
 & \leq  & \sum _{n=0}^{\infty }2^{-nq_{o}/2}\frac{\left|\xi \right|_{-p_{o}}^{n}}{n!}=\exp \left(\frac{\left|\xi \right|_{-p_{o}}}{2^{q_{o}/2}}\right)\, ,\label{eq:rho_{e}stimate}
\end{eqnarray}
so that $\rho _{\mu }(-\xi )\in \hmpqo $, and we can define the bilinear
form

\[
\check{B}_{\mu ,\nu }(\xi ,\eta ):=\bbform{\rho _{\mu }(-\xi )}{Be_{\nu }(\eta ;\cdot )}_{\mu }\]
 on the open cylinder $\Nc \times U_{p,q}\subset \Nc \times \Nc $. 

We call that bilinear form the \textbf{local $CS_{\mu \nu }$-symbol}
of the operator $B$.

\begin{example}
Let $B:\test \to \test $ be the continuous operator defined by ${B\bform{P_{n,\nu }}{\varphi _{n}}=\bform{\Pn }{\varphi _{n}}}$,
$\varphi _{n}\in \Nct{n}$, where $\mu ,\nu $ are two non-degenerate
measures with analytic Laplace transform. Then, $\check{B}_{\mu \nu }(\xi ,\eta )=\LL \rho _{\mu }(-\xi ;\cdot )|e_{\mu }(\eta ;\cdot )\GG _{\mu }=e^{\bform{\xi }{\eta }}$.
\end{example}
\begin{lem}
\label{lem:C-transform-hpq}Let $p_{o}$, and $q_{o}$ be non-negative
integers. Then for any ${\varphi \in \hpqo }$ the function ${C_{\mu }\varphi (\xi ):=\bbform{\rho _{\mu }(-\xi )}{\varphi }_{\mu }}$
is holomorphic on $H_{-p_{o},\mathbb{C}}$ and hence on $\Nc \subset H_{-p_{o},\mathbb{C}}$. 
\end{lem}
\begin{proof}
Let $\left\{ \varphi _{(n)}\in \test \right\} $ be a sequence of
test functions converging to $\varphi $ in $\hpqo $. Since $\abs{\bbform{\rho _{\mu }(-\xi ;\cdot )}{\varphi -\varphi _{(n)}}_{\mu }}\leq \exp \left(\frac{\left|\xi \right|_{-p_{o}}}{2^{q_{o}/2}}\right)\left\Vert \varphi -\varphi _{(n)}\right\Vert _{p_{o},q_{o},\mu }$,
the sequence of holomorphic functions $C_{\mu }\varphi _{(n)}(\xi )$
converges to $C_{\mu }\varphi (\xi )$ in the compact-open topology.
Therefore, $C_{\mu }\varphi (\xi )$ is holomorphic.
\end{proof}
\begin{thm}
\label{thm: CS symbol} For any $p_{o}\geq 0$, $\varepsilon >0$,
there exists $C>0$ and a $0$-neighborhood $U\subset \Nc $, such
that the local $CS_{\mu \nu }$-symbol $\check{B}_{\mu \nu }(\xi ,\eta )$
is holomorphic in the cylinder $\Nc \times U$, and there is an estimate
\[
\abs{\check{B}_{\mu \nu }(\xi ,\eta )}\leq Ce^{\varepsilon \abs{\xi }_{-p_{o}}}\, .\]

\end{thm}
\begin{proof}
Choose an integer $q_{o}\geq 0$, such that $2^{-q_{o}/2}\leq \varepsilon $.
By the Kernel Theorem, there exist $p\geq 0$, $q\geq 0$, such that
the operator $B:\hpqnu \rightarrow \hpqo $ is continuous and there
is a well defined bilinear form $\check{B}_{\mu \nu }(\xi ,\eta ):=\bbform{\rho _{\mu }(-\xi ;\cdot )}{Be_{\nu }(\eta ;\cdot )}_{\mu }$
on the cylinder $\Nc \times U_{p,q}$, where ${U_{p,q}=\{\eta \in \mathcal{N}_{\mathbb{C}};\, \left|\eta \right|_{p}^{2}<2^{-q-1}\}}$. 

For fixed $\eta $, the function $\check{B}_{\mu \nu }(\xi ,\eta )=\left[C_{\mu }Be_{\nu }(\eta ;\cdot )\right](\xi )$
is holomorphic on the whole $\Nc $ by Lemma \ref{lem:C-transform-hpq}.
For fixed $\xi $, the function $\check{B}_{\mu \nu }(\xi ,\cdot )$
is the local $S_{\nu }$-transform of the distribution $B^{*}\rho _{\mu }(-\xi ;\cdot )\in \hmpqnu $,
and therefore $\check{B}_{\mu \nu }(\xi ,\cdot )\in \holo $. Therefore,
the functions of one complex variable \[
\lambda \mapsto \check{B}_{\mu \nu }(\xi _{0}+\lambda \xi ,\cdot ),w\mapsto \check{B}_{\mu \nu }(\cdot ,\eta _{0}+w\eta )\, ;\, \, \, \lambda \in \mathbb{C};\eta _{0}\in U_{p,q};\xi _{0},\xi ,\eta \in \Nc \]
 are holomorphic at $0$. By the Hartogs' theorem (see the textbook
of Dineen\cite{dineen81}), the function of two complex variables\[
\lambda ,w\mapsto \check{B}_{\mu \nu }(\xi _{0}+\lambda \xi ,\eta _{0}+w\eta )\, ;\, \, \, \lambda ,w\in \mathbb{C};\eta _{0}\in U_{p,q};\eta ,\xi _{0},\xi \in \Nc \]
 is holomorphic at $0$. Another words, $\check{B}_{\mu ,\nu }(\xi ,\eta )$
is G-holomorphic on $\Nc \times U_{p,q}$. 

Observe, that $\check{B}_{\mu ,\nu }(\xi ,\eta )$ is also locally
bounded on $\Nc \times U_{p,q}$, as there exists $C>0$ such that
\begin{eqnarray*}
\abs{\bbform{\rho _{\mu }(-\xi ;\cdot )}{e_{\nu }(\eta ;\cdot )}_{\mu }} & \leq  & \left\Vert \rho _{\mu }(-\xi ;\cdot )\right\Vert _{-p_{o},-q_{o},\mu }C\left\Vert e_{\nu }(\eta ;\cdot )\right\Vert _{p,q,\nu }\\
 & \leq  & C\exp \left(\frac{\abs{\xi }_{-p_{o}}}{2^{q_{o}/2}}\right)\frac{1}{\sqrt{1-2^{q}\left|\eta \right|_{p}^{2}}}\, ,
\end{eqnarray*}
where $C$ is the norm of the operator $B:\hpq \to \hpqo $. Therefore,
${\check{B}_{\mu \nu }(\xi ,\eta )}$ is holomorphic on ${\Nc \times U_{p,q}}$.
\end{proof}

\subsection{Reconstruction of Operators from Their Local $CS_{\mu \nu }$-Symbols }

Now, let $F(\xi ,\eta )$ be a complex valued function defined and
G-holomorphic on an open cylinder ${\Nc \times U_{0}\subset \Nc \times \Nc }$,
and let $p_{o}$ be a non-negative integer.

\begin{thm}
Let $F(\xi ,\eta )$ be a complex valued function defined and G-holomorphic
on an open cylinder ${\Nc \times U_{0}\subset \Nc \times \Nc }$,
and suppose that for any $p_{o}\geq 0$, $\varepsilon >0$ there exists
$C>0$ and a $0$-neighborhood $U\subset U_{0}$ such that \begin{equation}
\abs{F(\xi ,\eta )}\leq Ce^{\varepsilon \abs{\xi }_{-p_{o}}}\, ,\, \, \, \, \, \xi \in \Nc ,\eta \in U.\label{eq:Thm-CS-growth}\end{equation}
 Then, there exists a unique continuous operator $B:(\mathcal{N})^{1}\rightarrow \test $,
such that its $CS_{\mu \nu }$-symbol ${\check{B}_{\mu \nu }(\xi ,\eta )}$
coincides with $F(\xi ,\eta )$ on some open cylinder ${\mathcal{N}_{\mathbb{C}}\times U\subset \mathcal{N}_{\mathbb{C}}\times \mathcal{N}_{\mathbb{C}}}$.
\end{thm}
\begin{proof}
Let $p_{o}\geq 0$, $q_{o}\geq 0$. Choose $\bar{p}>p_{o}$ such that
the embedding $I_{\bar{p},p_{o}}:H_{\bar{p}}\rightarrow H_{p_{o}}$
is a Hilbert-Schmidt operator. Choose $\varepsilon >0$ such that
${(\varepsilon e)^{2}2^{q_{o}}\left\Vert I_{\bar{p},p_{o}}\right\Vert _{HS}^{2}<1}$.
By assumption, there is $C>0$ and a neighborhood $U\subset U_{0}$
such that ${\abs{F(\xi ,\eta )}\leq Ce^{\varepsilon \abs{\xi }_{-p_{o}}}}$
on $\mathcal{N}_{\mathbb{C}}\times U$. Therefore, the function $F$
is holomorphic on $\mathcal{N}_{\mathbb{C}}\times U$. 

Let $\delta >0$ be such that ${\{\eta \in \mathcal{N}_{\mathbb{C}};\left|\eta \right|_{\bar{p}}\leq \delta \}\subset U}$
.  For every $(\xi ,\eta )\in \mathcal{N}_{\mathbb{C}}\times \mathcal{N}_{\mathbb{C}}$,
$\left|\eta \right|_{\bar{p}}\leq 1$ the function of two complex
variables ${(s,t)\mapsto F(s\xi ,t\eta )}$ can be expanded into the
power series\[
F(s\xi ,t\eta )=\sum _{m,n=0}^{\infty }\frac{F^{(m,n)}(\xi ,..,\xi ;\eta ,...,\eta )}{m!n!}s^{m}t^{n}\]
 converging on the cylinder $\left\{ (s,t)\in \mathbb{C}^{2};\left|t\right|\leq \delta \right\} $.
Here $F^{(m,n)}$ is an $(m+n)$-linear form \begin{eqnarray*}
F^{(m,n)}(\xi _{1},..,\xi _{m};\eta _{1},...,\eta _{n}) & = & \frac{\partial }{\partial s_{1}}\cdots \frac{\partial }{\partial s_{m}}\frac{\partial }{\partial t_{1}}\cdots \frac{\partial }{\partial t_{n}}F(s_{1}\xi _{1}+..+s_{n}\xi _{m};t_{1}\eta _{1}+...t_{n}\eta _{n})
\end{eqnarray*}
 symmetric with respect to the first $m$ and the last $n$ variables.
For any $R>0$ we have the Cauchy formula\[
\frac{1}{m!n!}F^{(m,n)}(\xi ,..,\xi ;\eta ,...,\eta )=(2\pi i)^{-2}\int \limits _{\left|s\right|=R,\left|t\right|=\delta }\frac{F(s\xi ,t\eta )}{s^{m+1}t^{n+1}}\, ds\, dt\]
 which gives the following Cauchy inequality for $\left|\xi \right|_{-\bar{p}}\leq 1,\left|\eta \right|_{\bar{p}}\leq 1$
\begin{equation}
\left|\frac{1}{m!n!}F^{(m,n)}(\xi ,..,\xi ;\eta ,...,\eta )\right|\leq Ce^{\varepsilon R}R^{-m}\delta ^{-n}\, ,\label{eq:Cauchy-CS}\end{equation}
so that ${\sup \left\{ \left|\frac{1}{m!n!}F^{(m,n)}(\xi ,..,\xi ;\eta ,...,\eta )\right|;\abs{\xi }_{-\bar{p}}\leq 1,\abs{\eta }_{\bar{p}}\leq 1\right\} \leq Ce^{\varepsilon R}R^{-m}\delta ^{-n}}$.
Choose $R=\frac{m}{\varepsilon }$. Using the polarization identity
and inequality ${\frac{n^{n}}{n!}\leq e^{n}}$, we obtain the estimate\begin{eqnarray*}
\abs{\frac{1}{m!n!}F^{(m,n)}(\xi _{1},..,\xi _{m};\eta _{1},...,\eta _{n})} & \leq  & \frac{m^{m}}{m!}\frac{n^{n}}{n!}Ce^{\varepsilon \frac{m}{\varepsilon }}\left(\frac{m}{\varepsilon }\right)^{-m}\delta ^{-n}\\
 & \leq  & \frac{C}{m!}e^{n}(\varepsilon e)^{m}\delta ^{-n}
\end{eqnarray*}
 on the {}``unit polydisc'' ${\left\{ (\xi _{1},..,\xi _{m};\eta _{1},...,\eta _{n});\xi _{i}\in \mathcal{N}_{\mathbb{C}},\eta _{j}\in \mathcal{N}_{\mathbb{C}},\abs{\xi _{i}}_{-\bar{p}}\leq 1,\abs{\eta _{j}}_{\bar{p}}\leq 1\right\} }$.
For any ${\xi _{1},..,\xi _{m},\, \eta _{1},...,\eta _{n}\in \mathcal{N}_{\mathbb{C}}}$
we have the estimate\[
\abs{\frac{1}{m!n!}F^{(m,n)}(\xi _{1},..,\xi _{m};\eta _{1},...,\eta _{n})}\leq C\frac{(\varepsilon e)^{m}}{m!}\left(\frac{e}{\delta }\right)^{n}\prod _{i=1}^{m}\left|\xi _{i}\right|_{-\bar{p}}\prod _{j=1}^{n}\left|\theta _{j}\right|_{\bar{p}}\, .\]

By the kernel theorem, for all $p>\bar{p}$ such that the embedding
${I_{p,\bar{p}}:H_{p}\rightarrow H_{\bar{p}}}$ is a Hilbert-Schmidt
operators, there exist unique kernels $f_{m,n}\in H_{p_{o},\mathbf{C}}^{\hat{\otimes }m}\otimes H_{-p,\mathbf{C}}^{\hat{\otimes }n}$
such that \[
\frac{1}{m!n!}F^{(m,n)}(\xi _{1},..,\xi _{m};\eta _{1},...,\eta _{n})=\langle f_{m,n}|\xi _{1}\hat{\otimes }..\hat{\otimes }\xi _{m}\otimes \eta _{1}\hat{\otimes }...\hat{\otimes }\eta _{n}\rangle \, .\]

Moreover, we have the following norm estimate

\[
\left|f_{m,n}\right|_{p_{o},-p}\leq C\frac{\left(\varepsilon e\norm{I_{\bar{p},p_{o}}}_{HS}\right)^{m}}{m!}\left(\frac{e}{\delta }\norm{I_{p,\bar{p}}}_{HS}\right)^{n}\, .\]

For any $\phi _{n}\in H_{p,\mathbb{C}}^{\hat{\otimes }n}$ define
$\langle f_{m,n}|n!\phi _{n}\rangle \in H_{p_{o},\mathbb{C}}^{\hat{\otimes }m}$
by the formula ${\langle \langle f_{m,n}|n!\phi _{n}\rangle \bigm |\psi _{m}\rangle :=n!\langle f_{m,n}|\psi _{m}\otimes \phi _{n}\rangle }$.
The definition is correct as we have the following norm estimate \[
\left|\langle f_{m,n}|n!\phi _{n}\rangle \right|_{p_{o}}\leq n!C\frac{\left(\varepsilon e\norm{I_{\bar{p},p_{o}}}_{HS}\right)^{m}}{m!}\left(\frac{e}{\delta }\norm{I_{p,\bar{p}}}_{HS}\right)^{n}\left|\phi _{n}\right|_{p}\, ,\]
 which shows that $\phi _{n}\mapsto \langle f_{mn}|n!\phi _{n}\rangle $
defines a continuous linear operator ${H_{p,\mathbb{C}}^{\hat{\otimes }n}\to H_{p_{o},\mathbb{C}}^{\hat{\otimes }m}}$. 

For an element $\phi (\cdot )=\sum _{n=0}^{\infty }\langle \Pn (\cdot );\phi _{n}\rangle \in (H_{p,q,\nu })$,
where $\phi _{n}\in H_{p,\mathbb{C}}^{\hat{\otimes }n}$, put\[
b_{m}\phi :=\sum _{n=0}^{\infty }\langle f_{m,n}|n!\phi _{n}\rangle \, ,\, \, \, \, \, m=0,1\ldots .\]

Using the Schwartz inequality we obtain the following estimate\begin{eqnarray}
\abs{b_{m}\phi }_{p_{o}}^{2} & \leq  & \left(\sum _{n=0}^{\infty }\left|\langle f_{m,n}|n!\phi _{n}\rangle \right|_{p_{o}}\right)^{2}\nonumber \\
 & \leq  & \left(\sum _{n=0}^{\infty }n!\left|f_{m,n}\right|_{p_{o},-p}\left|\phi _{n}\right|_{p}\right)^{2}\nonumber \\
 & = & \left(\sum _{n=0}^{\infty }n!2^{nq/2}\left|\phi _{n}\right|_{p}2^{-nq/2}\left|f_{m,n}\right|_{p_{o},-p}\right)^{2}\nonumber \\
 & \leq  & \sum _{n=0}^{\infty }(n!)^{2}2^{nq}\left|\phi _{n}\right|_{p}^{2}\sum _{n=0}^{\infty }2^{-nq}\left|f_{m,n}\right|_{p_{o},-p}^{2}\nonumber \\
 & \leq  & \left\Vert \phi \right\Vert _{p,q,\nu }^{2}C\frac{\left(\varepsilon e\norm{I_{\bar{p},p_{o}}}_{HS}\right)^{2m}}{(m!)^{2}}\sum _{n=0}^{\infty }2^{-nq}\left(\frac{e}{\delta }\norm{I_{p,\bar{p}}}_{HS}\right)^{2n}\nonumber \\
 & = & \left\Vert \phi \right\Vert _{p,q,\nu }^{2}C\frac{\left(\varepsilon e\norm{I_{\bar{p},p_{o}}}_{HS}\right)^{2m}}{(m!)^{2}}\left(1-2^{-q}\left(\frac{e}{\delta }\right)^{2}\left\Vert I_{p,\bar{p}}\right\Vert _{HS}^{2}\right)^{-1}\, ,\label{eq:cs-nsum-bmn}
\end{eqnarray}
 where $q$ is any positive integer such that $2^{-q}\left(\frac{e}{\delta }\right)^{2}\left\Vert I_{p,\bar{p}}\right\Vert _{HS}^{2}<1$.
That estimate shows that $b_{m}\phi \in H_{p_{o},\mathbb{C}}^{\hat{\otimes }m}$.
It also follows that\begin{eqnarray*}
 &  & \sum _{m=0}^{\infty }(m!)^{2}2^{mq_{o}}\left|b_{m}\phi \right|_{p_{o}}^{2}\\
 & \leq  & C\left\Vert \phi \right\Vert _{p,q,\nu }^{2}\left(1-(\varepsilon e)^{2}2^{q_{o}}\left\Vert I_{\bar{p},p_{o}}\right\Vert _{HS}^{2}\right)^{-1}\left(1-2^{-q}\left(\frac{e}{\delta }\right)^{2}\left\Vert I_{p,\bar{p}}\right\Vert _{HS}^{2}\right)^{-1}\, .
\end{eqnarray*}
 Thus, for any $p_{o},q_{o}\geq 0$ there exist $p\geq p_{o},q\geq 0$
such that\begin{eqnarray*}
B:\phi \mapsto \sum _{m=0}^{\infty }P_{m,\mu }(b_{m}\phi ) &  & 
\end{eqnarray*}
 is a continuous linear operator from $\hpqnu $ to $(H_{p_{o},q_{o},\mu })$.
This operator is unique, as the set $\left\{ e_{\nu }(\eta ,\cdot );\eta \in \Nc ,\abs{\eta }_{p}\leq \delta \right\} $
is total in $\hpqnu $ by Lemma \ref{lem:total-set}, the set $\left\{ \rho _{\mu }(\xi ,\cdot );\xi \in \Nc \right\} $
is total in $(H_{-p_{o},-q_{o},\mu })$ by Theorem \ref{thm:c-transfrom},
and there is the following identity for the $CS_{\mu \nu }$ -symbol\begin{eqnarray*}
\check{B}_{\mu \nu }(\xi ,\eta ) & = & \LL \rho _{\mu }(-\xi ;\cdot )|Be_{\nu }(\eta ;\cdot )\GG _{\mu }\\
 & = & \sum _{m=0}^{\infty }m!\left\langle \left(\sum _{n=0}^{\infty }n!\bform{f_{m,n}}{\frac{\xi ^{\otimes n}}{n!}}\right)\Biggl |\Biggr .\frac{\eta ^{\otimes m}}{m!}\right\rangle \\
 & = & \sum _{m,n=0}^{\infty }\langle f_{m,n}|\xi ^{\otimes m}\otimes \eta ^{\otimes n}\rangle \\
 & = & F(\xi ,\eta )\, ,\, \, \, \, \, \textrm{for}\, (\xi ,\eta )\in \mathcal{N}_{\mathbb{C}}\times \mathcal{N}_{\mathbb{C}}\, ,\, \left|\eta \right|_{p}\leq \left|\eta \right|_{\bar{p}}<\delta \, .
\end{eqnarray*}

It follows that $B\test \subset \test $, and $B:\test \to \test $
is the unique continuous operator such that its $CS_{\mu \nu }$-symbol
${\check{B}_{\mu \nu }(\xi ,\eta )}$ coincides with $F(\xi ,\eta )$
on some open cylinder ${\mathcal{N}_{\mathbb{C}}\times U\subset \mathcal{N}_{\mathbb{C}}\times \mathcal{N}_{\mathbb{C}}}$.
\end{proof}

\section*{Acknowledgement}

I would like to express my gratitude to Prof. Alexander Dynin for
his encouragement and helpful discussions. 

\bibliographystyle{abbrv}
\bibliography{liter}

\end{document}